\documentclass[winedt,
                       yap
                            ]{pit}
\usepackage{srcltx}
\usepackage{amsmath}
\usepackage{amssymb}
\usepackage{array}
\usepackage{longtable}
\usepackage{tabularx}
\usepackage{epsfig}
\usepackage{psfrag}
\usepackage{enumerate}
\usepackage{url}
\usepackage{mathrsfs}
\usepackage{dsfont}     % instead of bbm    !!!

\usepackage{euscript}
\usepackage{multirow}
\let\kappa\varkappa

\let\epsilon\varepsilon

\let\phi\varphi
\let\upn=\textup

\def\etal{{\noexpand\lowercase{\noexpand\normalfont et al.}}}
\newcolumntype{C}{>{\centering\arraybackslash}X}
\def\E{\mathop{\mbox{\bf{E}}}\nolimits}
\def\P{\mathop{\mbox{\bf{P}}}\nolimits}

\def\mmid{\mathchoice{\hskip1.5pt|\hskip1.5pt}{\hskip1.5pt|\hskip1.5pt}{\hskip0.5pt|\hskip0.5pt}{\hskip0.3pt|\hskip0.3pt}}

\makeatletter
\renewcommand{\pod}[1]{\if@display\mkern10mu\else\mkern6mu\fi(#1)}
\renewcommand{\pmod}[1]{\pod{{\operator@font mod}\mkern6mu#1}}
\makeatother
\renewcommand{\int}{\intop\limits}
\def\?{\mbox{\bf ??}\marginpar{\Large\bf ?}}

\allowdisplaybreaks

\begin{document}
\frontmatter          % for the preliminaries
%
%Here may be written anything what you want to be included BEFORE
%the table of contents and the main body of the journal issue
%
%
%
%Below variables should be specified for EACH issue of journal
%
\JournalName{Problems of Information Transmission}% by default -- Problems of Information Transmission
\JournalISSNCode{0032-9460}
\TransYearOfIssue{2017}
\TransCopyrightYear{2017}
\TransVolumeNo{53}
\TransIssueNo{4}
\OrigYearOfIssue{2017}
\OrigCopyrightYear{2017}
\OrigVolumeNo{53}
\OrigIssueNo{4}

%\renewcommand{\baselinestretch}{1.05}
%\tableofcontents % <--- use when preparing the final version
%
%\renewcommand{\baselinestretch}{1}
\mainmatter              % start of the contributions
\setcounter{page}{307}
%
%Below variables should be specified for EACH ARTICLE before the %\input command
%or inside the corresponding input article file
%
%\title{}% OBLIGATORY!
%\author{}% OBLIGATORY!
%\institute{}% OBLIGATORY!
%\received{}% OBLIGATORY!
%\titlerunning{}
%\authorrunning{}
%\toctitle{}
%\tocauthor{}
%\OrigPages{}% OBLIGATORY!
%\OrigCopyrightedAuthors{}% OBLIGATORY!

%\Rubrika{Information Theory}% The Rubrika name if necessary

\Rubrika{Coding Theory}

\Rubrika{Methods of Signal Processing}

\OrigPages{69--83}

\OrigCopyrightedAuthors{E.V.~Burnaev, G.K.~Golubev}

\title{On One Problem in Multichannel Signal Detection}

%\titlerunning{}
\authorrunning{Burnaev, Golubev}

\received{Received June 19, 2017; in final form, September 20, 2017}

\author{E. V.
Burnaev\inst{2}\dovesok{$^,$\kern1pt}\inst{1}\dovesok{$^,$}\thanks{Supported in part
by the Russian Foundation for Basic Research, project
no.~16-29-09649-ofi-m.}\dovesok{$^{,\kern1pt*}$} \and G. K.
Golubev\inst{1}\dovesok{$^,$\kern1pt}\inst{3}\dovesok{$^,$}\thanks{The research was
carried out at the Institute for Information Transmission Problems of the Russian
Academy of Sciences at the expense of the Russian Science Foundation, project
no.~14-50-00150.}\dovesok{$^{,\kern1pt**}$}}

\institute{Kharkevich Institute for Information Transmission Problems,\\ Russian
Academy of Sciences, Moscow, Russia \and Skolkovo Institute of Science and Technology
(Skoltech), Skolkovo, Russia \and CNRS, Aix-Marseille Universit\'e, I2M, UMR 7353,
Marseille, France\\ \email{$^*$e.burnaev@skoltech.ru}\upn, $^{**}${\tt
golubev.yuri@gmail.com}}

\doi{17040056}

\renewcommand{\aboveinst}{12pt}
\maketitle \setcounter{footnote}{2}
\renewcommand{\aboveinst}{6pt}

\begin{abstract}
We consider a statistical problem of  detection of a signal with unknown energy in a multi-channel system, observed in a Gaussian noise. We assume that the signal can appear in the $k$-th channel with a known small prior probability $\bar{\pi}_k$. Using noisy observations from all channels we would like to detect whether the signal is presented in one of the channels or we observe pure noise. In our work we describe and compare statistical properties of maximum posterior probability test and optimal Bayes test.
In particular, for these tests we obtain limiting distributions of test statistics and define sets of their non-detectable signals.
\end{abstract}

\section{Introduction}
We consider one of the basic problems of signal detection in a multi-channel system. From mathematical point of view we have to check a simple hypothesis $\mathrm{H}_0$, according to which the observed vector $Y\in\mathbb{R}^\infty$ can be represented as a discrete white noise
\begin{equation}\label{eq.1}
\mathrm{H}_0\colon\, Y=\sigma\xi,
\end{equation}
where $\xi=(\xi_i,\xi_2,\ldots)^\top$ is a standard Gaussian white noise, i.e. this is a vector in $\mathbb{R}^\infty$, $\sigma>0$ is  a known noise level.

Hypothesis alternative to $\mathrm{H}_0$ has the form
\begin{equation}\label{eq.2}
\mathrm{H}_1\colon\, Y=S+\sigma\xi,\quad S\in\mathbb{S},
\end{equation}
where $\mathbb{S}$ is a subset of signals in $\mathbb{R}^\infty$, which can have only one non-zero component. In other words, let $\mathbb{S}_k$ be a linear vector subspace in~$\mathbb{R}^\infty$, such that all coordinates except $k$-th are equal to zero. Then
\[
\mathbb{S}=\bigcup_{k=1}^\infty \mathbb{S}_k.
\]
We assume that the signal $S\in\mathbb{S}$ is random and independent of $\xi$, such that
\[
\P\bigl\{S\in\mathbb{S}_k \bigr\}=\bar{\pi}_k,
\]
where prior probabilities $\bar{\pi}_k$ are known.

Thus the problem can be formulated as follows: using observations $Y$ we want to check simple hypothesis $\mathrm{H}_0$ against complex alternative hypothesis $\mathrm{H}_1$.

Let us recall that a statistical test can be represented by any measurable function~$\phi(Y)$, taking values from the interval $[0,1]$. In the sequel for simplicity we assume that this function can take only two values $\{0,1\}$: if $\phi(Y)=0$, then we accept hypothesis $\mathrm{H}_0$, and if $\phi(Y)=1$, then we accept alternative hypothesis~$\mathrm{H}_1$.

Efficiency of any statistical test is measured by probabilities of errors of the first kind $\alpha_\phi$ (false alarm probability) and of the second kind $\beta_\phi$ (missing of a target), which can be defined as follows:
\[
\alpha_\phi =\P_0\{\phi(Y)=1\},\qquad \beta_\phi(S)=\P_S\{\phi(Y)=0\},
\]
where $\P_0$ is a probability measure of observations $Y$ from \eqref{eq.1}, and $\P_S$  is a probability measure of observations $Y$ from \eqref{eq.2} given fixed $S$.

Usually given a fixed false alarm rate, we would like to construct the test, which minimizes probability of missing a target signal. Unfortunately, we can not solve this problem in general case, since probability of the error of the second kind depends on $S$. However,  we can always  construct a statistical test, which for a given false alarm rate minimize average probability of the error of the second kind
\[
\bar{\beta}_{\phi}(Q)= \int_{\mathbb{S}} Q(S)\beta_\phi(S)\, dS,
\]
where positive function $Q(\cdot)$ is such that $\int_{\mathbb{S}} Q(S)\, dS=1.$ We would like to stress that probability density
 $Q(\cdot)$ should contain prior information about signal $S$.  Since we assume that $S\in\mathbb{S}_k$ with probability $\bar{\pi}_k$, and $\mathbb{S}_k$ is a one-dimensional subspace, then
\[
\int_{\mathbb{S}} Q(S)\beta_\phi(S)\, dS= \sum_{k=1}^\infty
\bar{\pi}_k\int_{-\infty}^\infty q_k(S_k)\beta_\phi(0,\ldots,0,S_k,0,\ldots)\, dS_k,
\]
where $q_k(\cdot)$ is a prior density of signal $S_k$ distribution in $k$-th channel.

From elementary course on mathematical statistics (Neyman-Pearson lemma) it is well-known that the test, minimizing average probability of the error of the second kind has the following form:
\begin{equation}\label{eq.3}
\phi(Y)=\mathbf{1}\Biggl\{\sum_{k=1}^\infty \pi_k \int_{-\infty}^\infty q_k(s)
l(s;Y_k)\, ds\ge t_\alpha \Biggr\},
\end{equation}
where
\[
l(s;Y_k)= \exp\biggl(-\frac{s^2-2sY_k}{2\sigma^2}\biggr)
\]
is a likelihood ratio for $k$-th channel, and critical level $t_\alpha$ is set in order to provide given false alarm probability $\alpha$.

Since as a rule we do not have any  prior information about distribution of non-zero component of the signal, then mathematically we can represent this fact e.g. by assuming that this component has a Gaussian distribution with a big variance, i.e.
\[
q_k(s)=\frac{1}{\sqrt{2\pi}A}\exp\biggl(-\frac{s^2}{2A^2}\biggr),
\]
where $A\gg \sigma$. Then, integrating in \eqref{eq.3}, we get that
\[
\int_{-\infty}^\infty
q_k(s)l(s;Y_k)ds=\frac{1}{\sqrt{1+A^2/\sigma^2}}\exp\biggl[\frac{Y_k^2}{2\sigma^2(1
+\sigma^2/ A^2)}\biggr],
\]
and for $A^2/\sigma^2\to\infty$ we get the following Bayes test:
\[%\label{s1.e2}
\phi^\circ(Y)=\mathbf{1}\Biggl\{\sum_{k=1}^\infty\bar{\pi}_k
\exp\biggl(\frac{Y_k^2}{2\sigma^2}\biggr)\ge t_\alpha^\circ(\bar{\pi})\Biggr\},
\]
where the critical value $t_\alpha^\circ(\bar{\pi})$ is set in order to provide given false alarm probability $\alpha$, or, in other words, the critical value is a solution of the equation
\begin{equation}\label{eq.4}
\P\Biggl\{\sum_{k=1}^\infty\bar{\pi}_k \exp\biggl(\frac{\xi_k^2}{2}\biggr) \ge
t_\alpha^\circ(\bar{\pi})\Biggr\}=\alpha,
\end{equation}
where here and elsewhere  $\bar{\pi} =(\pi_1,\pi_2,\ldots)^\top$ is a vector of prior probabilities. Let us note that strictly speaking this test is a Bayes test with improper prior distribution, but for brevity we will call it simply Bayes test.

In practice besides Bayes test we often use Maximum A Posteriori test (MAP)
\[
%\label{eq.5}
\phi^*(Y)=\mathbf{1}\biggl\{\max_{i\ge1}\biggl[\bar{\pi}_i
\exp\biggl(\frac{Y_i^2}{2\sigma^2}\biggr)\biggr]\ge t_\alpha^*(\bar{\pi}) \biggr\},
\]
where critical level $t_\alpha^*(\bar{\pi})$ is selected in order to provide given false alarm probability $\alpha$, or, in other words, the critical value is a solution of the equation
\begin{equation}\label{eq.6}
\P\biggl\{\max_{i\ge 1}\biggl[\bar{\pi}_i
\exp\biggl(\frac{\xi_i^2}{2}\biggr)\biggr]\ge t_\alpha^*(\bar{\pi})\biggr\}=\alpha.
\end{equation}

The aim of this work is to find out in which way and to what extent Bayes test is better than MAP test. For this we will use additional assumption about prior probabilities $\bar{\pi}_k$. We assume that
\begin{equation}\label{eq.7}
\bar{\pi}_k=\bar{\pi}_k^n=\bar{\pi}\Bigl(\frac{k}{n}\Bigr)\Bigm/ \sum_{s=1}^\infty
\bar{\pi}\Bigl(\frac{s}{n}\Bigr),
\end{equation}
where $\bar{\pi}(x)$, $x\in\mathbb{R}^+$, is a non-negative bounded function, such that
\begin{equation}\label{eq.8}
\int_0^\infty \bar{\pi}(x)\, dx =1,\qquad H(\bar{\pi})=\int_0^\infty \bar{\pi}(x)\log
\frac{1}{\bar{\pi}(x)}\, dx <\infty.
\end{equation}
In other words, this assumption means that prior probabilities are small, having the order $n^{-1}$, but at the same time the entropy of the prior distribution is bounded by $\log(n)+C$, where $C<\infty$ for any $n>1$. In fact, value $n$ is an effective dimension of the problem, and in the subsequent considerations we consider properties of statistical tests given that $n\to\infty$.

Problem of signal detection in multi-channel systems has numerous technical applications and rich history. Various statistical problem statements and formulations of this problem are considered e.g. in \cite{Ba}.

Detection of signal with known entropy in Gaussian channels using maximum likelihood approach is studied in details in 
\cite[Section~8.2]{G}.

It seems that one of the first mathematical works about Bayesian signal detection for multi-channel systems is \cite{D}, in which they studied statistical model, composed of $n$ Rayleigh channels. Problem of Bayesian signal detection with known entropy in Gaussian channels was considered in \cite{BB}. In this paper they assumed that signal can appear in one of  $n$ channels with equal prior probabilities.

Let us also note that the monograph \cite{IS} contains many interesting and useful facts about detection of signals in multi-channel systems with Gaussian noises.

Let us stress that in this paper we investigate a situation, when prior probabilities of a signal, observed in different channels, are different, and energy of the signal is unknown and is a nuisance parameter. Since statistical problem of signal detection in multi-channel system is a high-dimensional problem, then opposed to low-dimensional problems its solution significantly depends on available prior information about detectable signals, and so results, provided in this paper, differ significantly from known results of papers, listed above.

The work has the following structure. Basic statistical properties of the MAP test and the Bayes test are provided in Sections~2,~3. Proofs of theorems are provided in Appendix.

\section{MAP test}
We get the following result about the critical level of the MAP test (see \eqref{eq.6}).

\begin{theorem}\label{lemma1} For $n\to\infty$
\begin{equation}\label{eq.9}
\log[t_\alpha^*(\bar{\pi}^n)]=\log\frac{1}{\sqrt{\pi}\alpha}
-\frac{1}{2}\log\bigg[\log\biggl(\frac{n}{\sqrt{\pi}\alpha}\biggr)\biggr]+o(1).
\end{equation}
\end{theorem}

\begin{figure}[tp]
\centering\vskip3pt
\begin{psfrags}
\psfrag{0.5}{\hskip-2pt\scriptsize$0.5$} \psfrag{0.55}{\hskip-3pt\scriptsize$0.55$}
\psfrag{0.6}{\hskip-2pt\scriptsize$0.6$} \psfrag{0.65}{\hskip-3pt\scriptsize$0.65$}
\psfrag{0.7}{\hskip-2pt\scriptsize$0.7$} \psfrag{0.75}{\hskip-3pt\scriptsize$0.75$}
\psfrag{0.8}{\hskip-2pt\scriptsize$0.8$} \psfrag{0.85}{\hskip-3pt\scriptsize$0.85$}
\psfrag{0.9}{\hskip-2pt\scriptsize$0.9$} \psfrag{0.95}{\hskip-3pt\scriptsize$0.95$}
\psfrag{1}{\hskip-2pt\scriptsize$1$} \psfrag{-0.35}{\hskip-10pt\scriptsize$-0.35$}
\psfrag{-0.3}{\hskip-8.5pt\scriptsize$-0.3$}
\psfrag{-0.25}{\hskip-10pt\scriptsize$-0.25$}
\psfrag{-0.2}{\hskip-8.5pt\scriptsize$-0.2$}
\psfrag{-0.15}{\hskip-10pt\scriptsize$-0.15$}
\psfrag{-0.1}{\hskip-8.5pt\scriptsize$-0.1$}
\psfrag{-0.05}{\hskip-10pt\scriptsize$-0.05$} \psfrag{0}{\hskip-1.5pt\scriptsize$0$}
\psfrag{0.05}{\hskip-5.5pt\scriptsize$0.05$} \psfrag{0.1}{\hskip-4pt\scriptsize$0.1$}
\psfrag{0.15}{\hskip-5.5pt\scriptsize$0.15$} \psfrag{n=400}{\scriptsize$n{=}400$}
\psfrag{n=40}{\scriptsize$n{=}40$} \psfrag{1-}{\hskip-1.5ex$1-\alpha$}
\includegraphics[width=\textwidth,height=0.4\textheight]{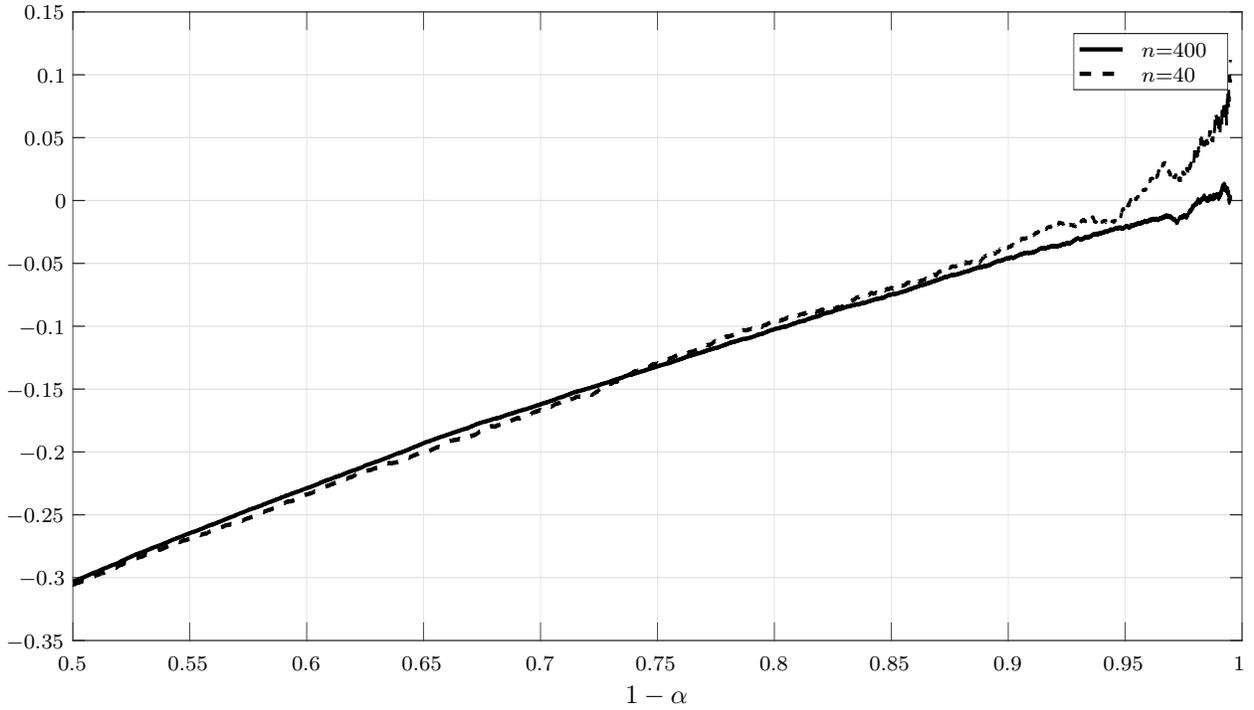}
\caption{ $\log[t_\alpha^*(\bar{\pi}^n)]$}\label{Fig.1}
\end{psfrags}
\end{figure}

\begin{singleremark}
Although convergence speed in \eqref{eq.9} is very low, still this formula is appropriate for applications. In Figure~\ref{Fig.1} we plot error of approximation
\[
\Delta(\alpha,n)=\log[t_\alpha^*(\bar{\pi}^n)]
-\log\frac{1}{\sqrt{\pi}\alpha}+\frac12 \log\biggl(\log
\frac{n}{\sqrt{\pi}\alpha}\biggr)
\]
as a function of $1-\alpha\in[0.5;0.995]$ for $n=40$ and $n=400$. We use uniform (on $[0,1]$) prior $\bar{\pi}(\cdot)$, and we estimate the critical value $t_\alpha^*(\bar{\pi}^n)$ by the Monte-Carlo method with $10^6$ random samples. Also let us note that usually in practice we are interested in small false alarm rates, i.e. $\alpha\le 0.05$. In Figure~\ref{Fig.1} we can see that for such values of false alarm rate error of approximation is small and decreases when $\alpha\to 0$.
\end{singleremark}

In order to describe statistical properties of the MAP test, related to the error of the second kind, we need to introduce some additional notations. Let us define the following parallelepiped in $\mathbb{R}^\infty$:
\begin{equation}\label{eq.10}
\Pi_{\bar{\pi}^n,\alpha}^*=\biggl\{x\in\mathbb{R}^\infty:\: x_i^2\le 2\sigma^2
\biggl[\log\frac{1}{\bar{\pi}_i^n}+\log\frac{1}{\sqrt{\pi}\alpha}
-\frac{1}{2}\log\biggl(\log \frac{n}{\sqrt{\pi}\alpha}\biggr)\biggr]\biggr\}.
\end{equation}

The following theorem essentially states that no signal from  $\Pi_{\bar{\pi}^n,\alpha}^* \cap \mathbb{S}$ can be detected by the MAP test.

\begin{theorem}\label{theorem1}
\upn The error of the second kind for the MAP test fulfills the following inequality
\[%\label{map}
\lim_{n\to\infty}\,\inf_{S\in\Pi_{\bar{\pi}^n,\alpha}^*\cap
\mathbb{S}}\,\beta_{\phi^*}(S)\ge\frac{1-\alpha}{2}.
\]
\end{theorem}

Let us note that due to this theorem average energy of a signal, which can not be detected by the MAP test, is equal to
\[
\sum_{k=1}^\infty \bar{\pi}_k^nS_k^2=2\sigma^2 \bigl\{\log[n
t_\alpha^*(\bar{\pi}^n)]+H(\bar{\pi}) +o(1) \bigr\},
\]
where the entropy $H(\bar{\pi})$ is defined in \eqref{eq.8}, and the critical value $t_\alpha^*(\bar{\pi}^n)$ is defined in \eqref{eq.9}.

\section{Bayes test}

In order to analyze the Bayes test first of all we need to investigate behavior of
\[
\Sigma_n(\xi)=\frac{1}{n}\sum_{i=1}^n\exp\biggl(\frac{\xi_i^2}{2}\biggr)
\]
when $n\to\infty$.

The following theorem, which, in essence, represents the main result of this paper, states that the distribution of $\Sigma_n(\xi)$ can be represented using the distribution of the following random variable
\begin{equation}\label{eq.11}
\zeta^\circ =\sum_{k=1}^\infty \biggl[
\biggl(\sum_{s=1}^ke_s\biggr)^{-1}-\frac{1}{k}\biggr]+\gamma.
\end{equation}
From now on $e_s$ are independent standard exponentially distributed random variables, $\gamma=0.577215\ldots\strut$ is an Euler constant.

Let us denote for brevity
\begin{equation}\label{eq.12}
b_n=\biggl[2\log \frac{n}{\sqrt{\pi \log(n)}}\biggr]^{1/2}.
\end{equation}

\begin{theorem}\label{main} For $n\to\infty$ we get that
\begin{equation}\label{eq.13}
\Sigma_n(\xi) \stackrel{\mathbf{P}}{=}\sqrt{\frac{2}{\pi}}\biggl(b_n
+\frac{\zeta^\circ}{b_n}\biggr) +o\biggl(\frac{1}{\sqrt{\log(n)}}\biggr).
\end{equation}
\end{theorem}

In \eqref{eq.13} and further in this paper for two sequences of random variables $\kappa_n$ and~$\kappa'_n$ 
notation
\[
\kappa_n\stackrel{\mathbf{P}}{=}\kappa'_n+o(r_n),\quad n\to \infty,
\]
means that there exist probability space on which these random variables are defined, and
\[
\lim_{n\to\infty}\P\biggl\{\frac{|\kappa_n-\kappa_n'|}{r_n}\ge\epsilon\biggr\}=0
\]
for any $\epsilon>0$.

\begin{figure}[tp]
\centering\vskip3pt
\begin{psfrags}
\psfrag{0.1}{\hskip-4pt\scriptsize$0.1$} \psfrag{0.2}{\hskip-4pt\scriptsize$0.2$}
\psfrag{0.3}{\hskip-4pt\scriptsize$0.3$} \psfrag{0.4}{\hskip-4pt\scriptsize$0.4$}
\psfrag{0.5}{\hskip-4pt\scriptsize$0.5$} \psfrag{0.6}{\hskip-4pt\scriptsize$0.6$}
\psfrag{0.7}{\hskip-4pt\scriptsize$0.7$} \psfrag{0.8}{\hskip-4pt\scriptsize$0.8$}
\psfrag{0.9}{\hskip-4pt\scriptsize$0.9$} \psfrag{1}{\hskip-2pt\scriptsize$1$}
\psfrag{-5}{\hskip-1pt\scriptsize$-5$} \psfrag{0}{\hskip-1.5pt\scriptsize$0$}
\psfrag{5}{\hskip-.5pt\scriptsize$5$} \psfrag{10}{\hskip-1.5pt\scriptsize$10$}
\psfrag{15}{\hskip-1.5pt\scriptsize$15$} \psfrag{20}{\hskip-1.5pt\scriptsize$20$}
\psfrag{25}{\hskip-1.5pt\scriptsize$25$} \psfrag{30}{\hskip-1.5pt\scriptsize$30$}
\psfrag{35}{\hskip-1.5pt\scriptsize$35$} \psfrag{40}{\hskip-1.5pt\scriptsize$40$}
\psfrag{45}{\hskip-2.5pt\scriptsize$45$} \psfrag{Distribution
function}{\hskip-4.5ex\raisebox{1.5ex}{ }} \psfrag{Distribution of
z}{\scriptsize$\zeta^\circ$} \psfrag{Distribution of 1e}{\scriptsize$1/e_1$}
\psfrag{X}{\hskip-.9ex$X$}
\includegraphics[width=\textwidth,height=0.4\textheight]{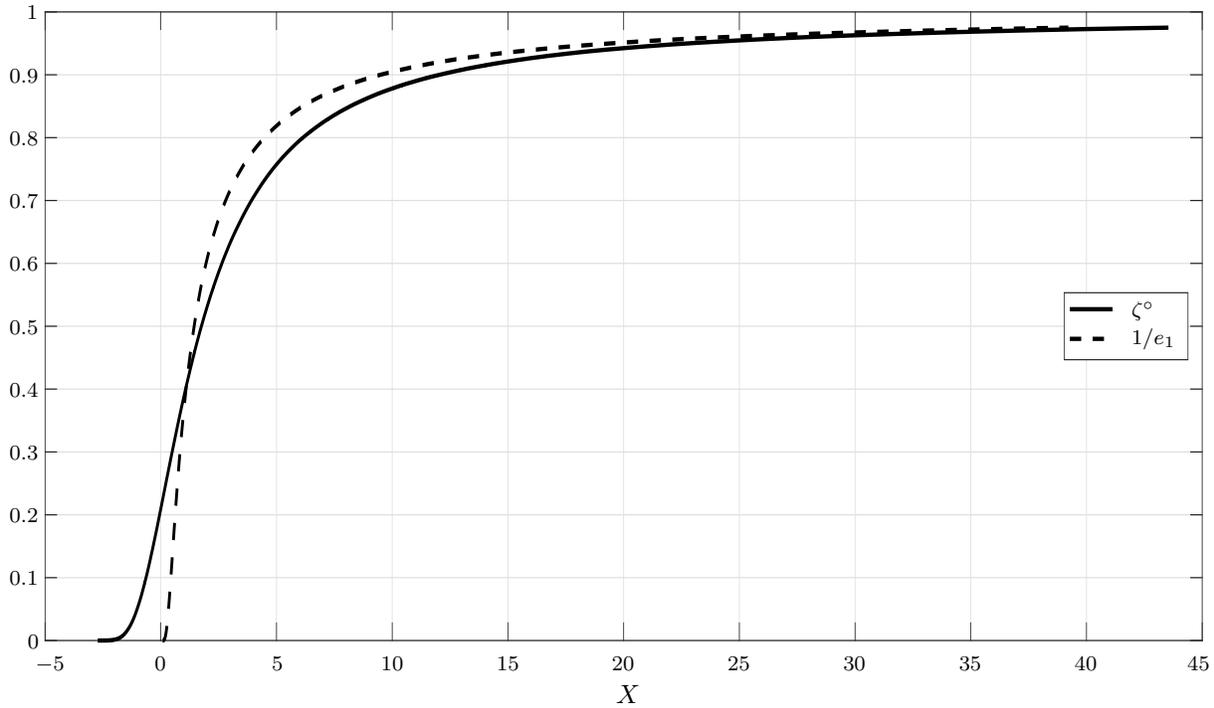}
\caption{ $\zeta^\circ$ \eqref{eq.11} ( ) $1/e_1$ ( )}\label{Fig.3}
\end{psfrags}
\end{figure}

In Figure.~\ref{Fig.3} we show distribution functions of random variables $\zeta^\circ$ and $1/e_1$. Let us note that 90\% of the mass of distribution of $\zeta^\circ$ is concentrated on the interval $[-1.02;26.01]$, namely
\[
\P\{\zeta^\circ>26.01\}=0.05\qquad \text{and}\qquad \P\{\zeta^\circ<-1.02\}=0.05.
\]
Besides that, the distribution of $\zeta^\circ$ has ``heavy tail''
\[
\P\{\zeta^\circ\ge x\}\asymp \frac{1}{x},\quad x\to \infty,
\]
which is clearly visible on the Figure.

Independent random variable, distributed as $\zeta^\circ$, have the following interesting property

\begin{theorem}\label{main2}
Let \/ $\zeta^\circ_i$\upn, $i=1,\ldots,p$\upn, be independent identically distributed random variables from \/ \eqref{eq.11}. Then for any \/
$\lambda_i>0$\upn, $i=1,\ldots,p$\upn, such that \/ $\sum\limits_{i=1}^p \lambda_i=1$\upn, we get that
\[
\P\Biggl\{\sum_{i=1}^p\lambda_i\zeta_i^\circ+\sum_{i=1}^p \lambda_i\log(\lambda_i)\le
x\Biggr\}= \P\{\zeta^\circ\le x\}.
\]
\end{theorem}
This result follows almost directly from theorem \ref{main}, so we omit its proof here.

From theorems \ref{main}, \ref{main2} the following interesting fact follows, describing distribution of statistics of the Bayes test under the null hypothesis.

\begin{theorem}\label{main3}
For $n\to\infty$
\[
\sum_{i=1}^\infty\bar{\pi}_i^n\exp\biggl(\frac{\xi_i^2}{2}\biggr)\stackrel{\mathbf{P}}{=}
\sqrt{\frac{2}{\pi}}\biggl[ b_n +\frac{\zeta^\circ+H(\bar{\pi})}{b_n}\biggr]
+o\biggl(\frac{1}{\sqrt{\log(n)}}\biggr);
\]
here \/ $\bar{\pi}_i^n$ and $H(\bar{\pi})$ are defined in \/ \eqref{eq.7} and \/ \eqref{eq.8} correspondingly.
\end{theorem}

Using this theorem we can get approximation of the critical value 
$t_\alpha^\circ(\bar{\pi}^n)$ from \eqref{eq.4}. Let us define quantile $t_\alpha^\circ$ as a solution of the equation
\begin{equation}\label{eq.14}
\P\{\zeta^\circ\ge t_\alpha^\circ\}=\alpha.
\end{equation}

Then
\begin{equation}\label{eq.15}
t_\alpha^\circ(\bar{\pi}^n)=\sqrt{\frac{2}{\pi}}\biggl[ b_n
+\frac{t_\alpha^\circ+H(\bar{\pi})}{b_n}\biggr]
+o\biggl(\frac{1}{\sqrt{\log(n)}}\biggr).
\end{equation}

\begin{figure}[tp]
\centering\vskip2pt
\begin{psfrags}
\psfrag{0.8}{\hskip-2pt\scriptsize$0.8$} \psfrag{0.82}{\hskip-3pt\scriptsize$0.82$}
\psfrag{0.84}{\hskip-3pt\scriptsize$0.84$} \psfrag{0.86}{\hskip-3pt\scriptsize$0.86$}
\psfrag{0.88}{\hskip-3pt\scriptsize$0.88$} \psfrag{0.9}{\hskip-2pt\scriptsize$0.9$}
\psfrag{0.92}{\hskip-3pt\scriptsize$0.92$} \psfrag{0.94}{\hskip-3pt\scriptsize$0.94$}
\psfrag{0.96}{\hskip-3pt\scriptsize$0.96$} \psfrag{0.98}{\hskip-3pt\scriptsize$0.98$}
\psfrag{1}{\hskip-2pt\scriptsize$1$} \psfrag{-0.7}{\hskip-8.5pt\scriptsize$-0.7$}
\psfrag{-0.6}{\hskip-8.5pt\scriptsize$-0.6$}
\psfrag{-0.5}{\hskip-8.5pt\scriptsize$-0.5$}
\psfrag{-0.4}{\hskip-8.5pt\scriptsize$-0.4$}
\psfrag{-0.3}{\hskip-8.5pt\scriptsize$-0.3$}
\psfrag{-0.2}{\hskip-8.5pt\scriptsize$-0.2$}
\psfrag{-0.1}{\hskip-8.5pt\scriptsize$-0.1$} \psfrag{0}{\hskip-1.5pt\scriptsize$0$}
\psfrag{n=400}{\scriptsize$n{=}400$} \psfrag{n=40}{\scriptsize$n{=}40$}
\psfrag{1-}{\hskip-1.5ex$1-\alpha$}
\includegraphics[width=\textwidth,height=0.4\textheight]{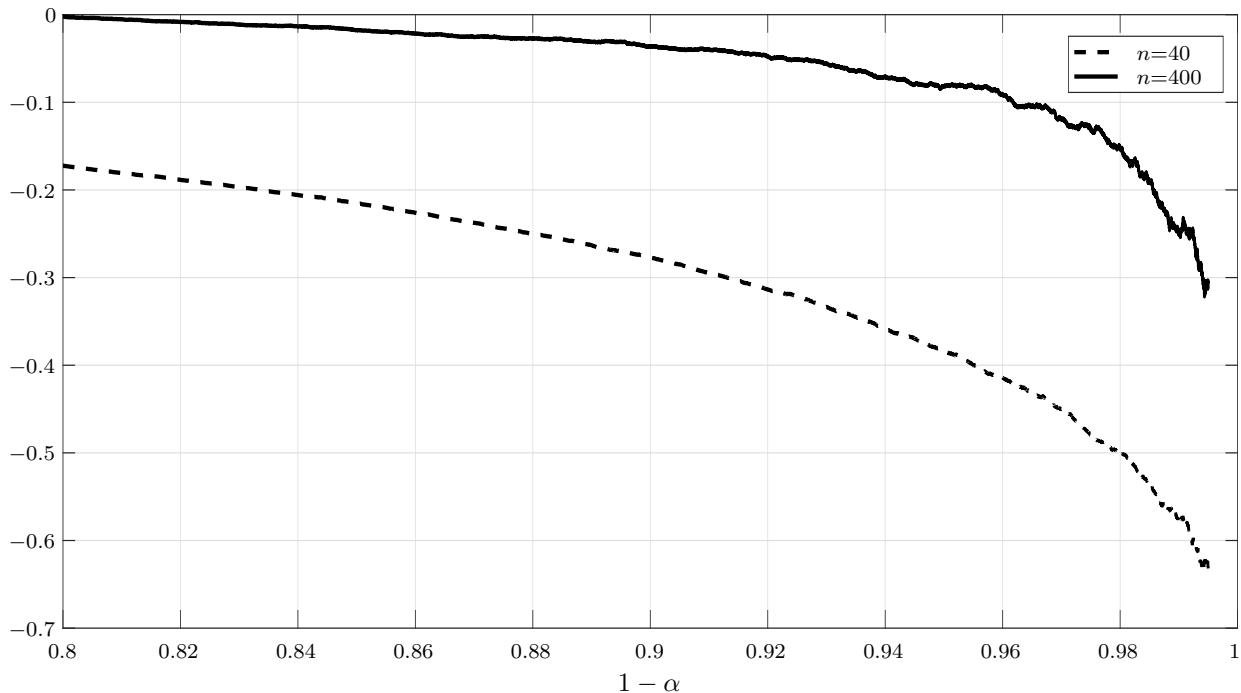}
\caption{ $t_\alpha^\circ(\bar{\pi}^n)$ \eqref{eq.15}}\label{Fig.2}
\end{psfrags}
\end{figure}

In Figure~\ref{Fig.2} we depict error of approximation of the critical level $t_\alpha^\circ(\bar{\pi}^n)$ using asymptotic expansion from \eqref{eq.15} for $\alpha\in [0.001;0.2]$ and
$n=40,400$. Despite the fact that convergence speed in \eqref{eq.15} is slow, we can see that this formula provides sufficiently accurate approximation.

As in the case with the MAP test, let us define the following parallelepiped:
\begin{equation}\label{eq.16}
\Pi_{\bar{\pi}^n}^\circ=\biggl\{x\in\mathbb{R}^\infty:\: x_i^2\le 2\sigma^2
\log\frac{1}{\bar{\pi}_i^n\sqrt{\pi\log(n)}}\biggr\}.
\end{equation}

\begin{theorem}\label{lemma-b}
\upn The error of the second kind for the Bayes test fulfills the following inequality
\[%\label{mmm}
\lim_{n\to\infty}\,\inf_{S\in\Pi_{\bar{\pi}^n}^\circ\cap
\mathbb{S}}\,\beta_{\phi^\circ}(S)\ge\frac{1-\alpha}{2}.
\]
\end{theorem}

Let us note that for big values of $n$ the difference between squares of the sides of parallelepipeds of non-detectable signals
 $\Pi_{\bar{\pi}^n,\alpha}^*$ and $\Pi_{\bar{\pi}^n}^\circ$ does not depend neither on side index nor on $n$ and is equal to (see
\eqref{eq.10} and \eqref{eq.16})
\[
\delta(\alpha)=2\sigma^2\log\frac{1}{\alpha}.
\]
The statistical meaning of this quantity is clear and means that the Bayes test in comparison with the MAP test can detect signals with  energies that are smaller by the value of $\delta(\alpha)$.

Another interesting feature of the Bayes test is that unlike the MAP test its parallelepiped of non-detectable signals does not depend on the probability of the false alarm $\alpha$.

\appendix{}

From this point on we will need simple probabilistic properties of independent random variables
$\exp(\xi_i^2/2)/n$, $i=1,\ldots,n$. Let us note that using integration by parts we can easily obtain that for $x\to\infty$
\begin{equation}\label{eq.17}
\P\biggl\{\frac{\xi_i^2}{2}\ge x \biggr\}=2 \P\Bigl\{\xi_i\ge\sqrt{2x}\Bigr\}
={(1+o(1))}\frac{{\rm e}^{-x}}{\sqrt{\pi x}}
\end{equation}
and therefore for $nx\to\infty$
\[
q_n(x)\stackrel{\rm def}{=}\P\biggl\{\frac1n\exp\biggl(\frac{\xi_i^2}{2}\biggr)\ge x
\biggr\}= \P\Bigl\{\frac{\xi_i^2}{2}\ge\log(nx)\Bigr\} =\frac{1+o(1)}{nx\sqrt{\pi
\log(nx)}}.
\]

Also, we can easily get asymptotic of the function $q^{-1}_n(x)$, which is inverse to $q_n(x)$:
\begin{equation}\label{eq.18}
q^{-1}_n(x)=\frac
{1+o(1)}{nx\sqrt{\pi}}\bigg[\log\biggl(\frac{1}{x\sqrt{\pi}}\biggr)\bigg]^{-1/2},
\quad x\to \infty.
\end{equation}

For brevity we denote by $x_{(i)}$, $i=1,\ldots,n$ values of $x_i$, $i=1,\ldots,n$  sorted in ascending order.

Let $U_i$, $i=1,\ldots,n$, be random variables, uniformly distributed on $[0,1]$. Thanks to Pyke theorem \cite{P}
\[
U_{(k)}\stackrel{\mathbf{P}}{=}\sum_{s=1}^k e_s \Bigm/ \sum_{s=1}^{n+1} e_s.
\]

Thus we can easily check that for $n\to\infty$
\begin{equation}\label{eq.19}
U_{(n-k)}\stackrel{\mathbf{P}}{=}
1-\frac{1}{n}\sum_{s=n-k+1}^{n+1}e_s+o\biggl(\frac{n-k}{n}\biggr).
\end{equation}

Then using \eqref{eq.19} and \eqref{eq.18}, we get that
\begin{multline}\label{eq.20}
\frac1n\exp\biggl(\frac{\xi^2_{(n-k)}}{2}\biggr)=q^{-1}_n\bigl(1-U_{(n-k)}\bigr)
=q^{-1}_n\Biggl(\frac{1}{n}\sum_{s=n-k+1}^{n+1}e_s\Biggr)\\
{\stackrel{\mathbf{P}}{=}}
(1+o(1))\Biggl(\sqrt{\pi}\sum_{s=n-k+1}^{n+1}e_s\Biggr)^{-1}
\Biggl[\log\Biggl(\frac{n}{\sqrt{\pi}}\Bigm/
\sum_{s=n-k+1}^{n+1}e_s\Biggr)\Biggr]^{-1/2}.\quad
\end{multline}

\proofof{Theorem~\ref{lemma1}} Thanks to \eqref{eq.17} we get that
\[
\begin{aligned}
\P\biggl\{\max_{i\ge 1}\biggl[\frac{\xi_i^2}{2}+\log(\bar{\pi}_i^n)\biggr]\ge
x\biggr\} &=1-\P\biggl\{\max_{i\ge
1}\biggl[\frac{\xi_i^2}{2}+\log(\bar{\pi}_i^n)\biggr]< x\biggr\}\\ &=
1-\prod_{i=1}^\infty\biggl[1- \P\biggl\{\frac{\xi_i^2}{2}\ge
x-\log(\bar{\pi}_i^n)\biggr\}\biggr]\\ &= 1-\prod_{i=1}^\infty
\biggl[1-(1+o(1))\bar{\pi}_i^n
\frac{\exp(-x)}{\sqrt{\pi}\sqrt{x-\log(\bar{\pi}_i^n)}}\biggr]\\ &=
1-\exp\Biggl\{\sum_{i=1}^\infty\log\bigg[1-(1+o(1))\bar{\pi}_i^n
\frac{\exp(-x)}{\sqrt{\pi}\sqrt{x-\log(\bar{\pi}_i)}}\biggr]\Biggr\}.
\end{aligned}
\]
Then using Taylor formula and \eqref{eq.7}, we continue this chain of equalities as follows:
\[
\begin{aligned}
\P\biggl\{\max_{i\ge 1}\biggl[\frac{\xi_i^2}{2}+\log(\bar{\pi}_i^n)\biggr]\ge
x\biggr\} &=1-\exp\biggl\{-(1+o(1))\frac{\exp(-x)}{\sqrt{\pi}}\sum_{i=1}^\infty
\frac{\bar{\pi}_i^n}{\sqrt{-\log(\bar{\pi}_i^n)+x}}\biggr\}\\ &= (1+o(1))
\frac{\exp(-x)}{\sqrt{\pi}}\sum_{i=1}^\infty
\frac{\bar{\pi}_i^n}{\sqrt{-\log(\bar{\pi}_i^n)+x}}\\ &=(1+o(1))
\frac{\exp(-x)}{\sqrt{\pi}}\int_0^\infty
\frac{\pi(t)}{\sqrt{\log(n)+x+\log[\pi(t)]}}\,dt\\
&=(1+o(1))\frac{\exp(-x)}{\sqrt{\pi[\log(n)+x]}}.
\end{aligned}
\]

From the last equality, solving the equation
\[
(1+o(1))\frac{\exp(-x)}{\sqrt{\pi[\log(n)+x]}}=\alpha,
\]
we obtain \eqref{eq.9}.\qed

\smallskip
\proofof{Theorem~\ref{theorem1}} Let us denote by $\nu$ an index of a channel, in which signal appears. Then in order to calculate error of the second kind for MAP test let us note that
\begin{multline}\label{eq.21}
\P\biggl\{\max_{k\ge 1}\biggl[\frac{Y_k^2}{2\sigma^2}+\log(\bar{\pi}_k^n)\biggr]
\le\log[t_\alpha^*(\bar{\pi}^n)] \,\Big|\,\nu=j\biggr\}\\ =\P\Biggl\{\max_{k\ne
j}\biggl(\frac{\xi_k^2}{2}+\log(\bar{\pi}_k^n)\biggr)\le\log[t_\alpha^*(\bar{\pi}^n)]
\cap\biggl[\frac{1}{2}\biggl(\frac
{S_j}{\sigma}+\xi_j\biggr)^2+\log(\bar{\pi}_j^n)\biggr]
\le\log[t_\alpha^*(\bar{\pi}^n)] \Biggr\}\\ =
\P\Biggl\{-S_j-\sigma\sqrt{2\log\frac{t_\alpha^*(\bar{\pi}^n)}{\bar{\pi}_j^n}}
\le\xi_j\le -S_j+\sigma
\sqrt{2\log\frac{t_\alpha^*(\bar{\pi}^n)}{\bar{\pi}_j^n}}\Biggl\}\\ \times
\P\biggl\{\max_{k\ne j}\biggl(\frac{\xi_k^2}{2}+\log(\bar{\pi}_k^n)\biggr)
\le\log[t_\alpha^*(\bar{\pi}^n)]\biggr\}.\quad
\end{multline}

Without loss of generality we can assume that $S_j >0$. Then it is obvious that
\begin{equation}\label{eq.22}
\lim_{n\to \infty}
\P\Biggl\{-S_j-\sigma\sqrt{2\log\frac{t_\alpha^*(\bar{\pi}^n)}{\bar{\pi}_j^n}}\le\xi_j\le
-S_j+\sigma
\sqrt{2\log\frac{t_\alpha^*(\bar{\pi}^n)}{\bar{\pi}_j^n}}\Biggl\}\ge\frac12.
\end{equation}

The last multiplier in the right part of \eqref{eq.21} can be bounded from below as follows. We will associate with the vector $\bar{\pi}^n$, belonging to a simplex in $\mathbb{R}^\infty$, the vector $\bar{\pi}^{(-j)}$ with coordinates, calculated as follows:
\begin{equation}\label{eq.23}
\bar{\pi}_k^{(-j)}=
\begin{cases}
\bar{\pi}_k^n/M^{(-j)}(\bar{\pi}^n), & k<s,\\
\bar{\pi}_{k+1}^n/M^{(-j)}(\bar{\pi}^n), & k\ge s,
\end{cases}
\end{equation}
where
\[
M^{(-j)}(\bar{\pi}^n)=\sum_{k\ne j}\bar{\pi}_k^n.
\]
In other words, we delete $j$-th coordinate from the vector $\bar{\pi}^n$  and then normalize this new vector in order it belongs to a simplex.
Since $\bar{\pi}_j^n=O(n^{-1})$, then from theorem \ref{lemma1} we get that for $n\to\infty$
\[
\log[t_\alpha^*(\bar{\pi}^n)]=\log[t_\alpha^*(\bar{\pi}^{(-j)})]+o(1),
\]
and therefore
\begin{multline*}
\lim_{n\to\infty}\P\biggl\{\max_{k\ne j}
\biggl(\frac{\xi_k^2}{2}+\log(\bar{\pi}_k^n)\biggr)
\le\log[t_\alpha^*(\bar{\pi}^n)]\biggr\}\\ =\lim_{n\to\infty} \P\biggl\{\max_{k\ne
j}\biggl(\frac{\xi_k^2}{2}+\log(\bar{\pi}_k^n)\biggr)
\le\log[t_\alpha^*(\bar{\pi}^{(-j)})]+o(1)\biggr\}=1-\alpha.
\end{multline*}
Thanks to this inequality together with
\eqref{eq.21} and \eqref{eq.22} we complete the proof.\qed

\smallskip
\proofof{Theorem~\ref{main}} Let us divide the sum $\Sigma_n(\xi)$ into two parts:
\begin{equation}\label{eq.24}
\Sigma_n(\xi)=\Sigma^{\rm d}_n(\xi)+\Sigma^{\rm r}_n(\xi),
\end{equation}
where
\[
\Sigma^{\rm d}_n(\xi)=\frac1n\sum_{i=1}^n\exp\biggl(\frac{\xi_i^2}{2}\biggr)
\mathbf{1}\bigl\{|\xi_i|< h_n\bigr\},\qquad \Sigma^{\rm r}_n(\xi)=\frac1n
\sum_{i=1}^n\exp\biggl(\frac{\xi_i^2}{2}\biggr) \mathbf{1}\bigl\{|\xi_i|\ge
h_n\bigr\}.
\]
We define the threshold $h_n$ as a root of the equation
\begin{equation}\label{eq.25}
\frac{1}{h_n}\exp\biggl(\frac{h_n^2}{2}\biggr)=\frac{n}{\sqrt{2\pi} M_n \log(n)};
\end{equation}
Here $M_n\to\infty$ for $n\to\infty$, but slower, namely
\begin{equation}\label{eq.26}
\lim_{n\to\infty}\frac{\log(M_n)\log[\log(n)]}{\log(n)}=0.
\end{equation}

The next theorem states that the distribution of the random variable $\Sigma^{\rm d}_n(\xi)$ degenerates for big $n$.

\begin{singlelemma}
For $n\to\infty$
\begin{equation}\label{eq.27}
\Sigma^{\rm d}_n(\xi)\stackrel{\mathbf{P}}{=}\sqrt{\frac{2}{\pi}}h_n
+o\biggl(\frac{1}{\sqrt{\log(n)}}\biggr).
\end{equation}
\end{singlelemma}

\begin{proof}
Let us note that
\[
\E \exp\biggl(\frac{\xi_i^2}{2}\biggr) \mathbf{1}\bigl\{|\xi_i|< h_n\bigr\}=
\frac{2}{\sqrt{2\pi}}\int_0^{h_n} \, dx=\frac{2h_n}{\sqrt{2\pi}},
\]
and so due to \eqref{eq.25}
\[
\begin{aligned}
\E \exp\bigl({\xi_i^2}\bigr)\mathbf{1}\bigl\{|\xi_i|< h_n\bigr\}&=
\frac{2}{\sqrt{2\pi}}\int_0^{h_n} \exp\biggl(\frac{x^2}{2}\biggr)\,dx\\ &=
(1+o(1))\frac{\sqrt{2}}{\sqrt{\pi} h_n} \exp\biggl(\frac{h_n^2}{2}\biggr)=
O\biggl(\frac{n}{M_n\log(n)}\biggr).
\end{aligned}
\]
Equality \eqref{eq.27} is an obvious consequence of these relations and \eqref{eq.26}.\qed
\end{proof}

Then we use a simple formula, which can be derived explicitly from ~\eqref{eq.25} and Taylor formula:
\begin{equation}\label{eq.28}
h_n=b_n-\frac{\log(M_n)}{b_n}+o\biggl(\frac{1}{\sqrt{\log(n)}}\biggr),
\end{equation}
where $b_n$ is defined in \eqref{eq.12}.

Let us consider random variable $\Sigma^{\rm r}_n(\xi)$, denote for the sake of brevity
\[
\mu_i=\frac{1}{n}\exp\biggl(\frac{\xi_i^2}{2}\biggr)
\]
and use simple relations (see \eqref{eq.25} and \eqref{eq.26})
\begin{equation}\label{eq.29}
\begin{aligned}[b]
\Sigma^{\rm r}_n(\xi)&=\sum_{k=0}^{n-1}\mu_{(n-k)}
\mathbf{1}\biggl\{\mu_{(n-k)}\ge\frac1n \exp\biggl(\frac{h_n^2}{2}\biggr)\biggr\}\\
&=\sum_{k=0}^{n-1}\mu_{(n-k)}\mathbf{1}\biggl\{\mu_{(n-k)} \ge\frac{h_n}{\sqrt{2\pi}
M_n \log(n)}\biggr\}\\ &=\sum_{k=0}^{n-1}\mu_{(n-k)}
\mathbf{1}\biggl\{\mu_{(n-k)}\ge\frac{1+o(1)}{\sqrt{\pi\log(n)} M_n}\biggr\}.
\end{aligned}
\end{equation}

Let us also define
\[
\mathcal{E}(i)=\sum_{k=1}^i e_k.
\]
Then from \eqref{eq.20} it follows that
\[
\mu_{(n-i)}\stackrel{\mathbf{P}}{=}\frac{1+o(1)}{\mathcal{E}(i)
\sqrt{\pi\log[n/\mathcal{E}(i)]}}.
\]
Let us define a ``stopping moment''
\[
\tau=\max \biggl\{k:\: \mu_{(n-k)}\ge\frac{1+o(1)}{\sqrt{\pi \log(n)} M_n}\biggr\}
=\max \Biggl\{k:\: \mathcal{E}(k)
\sqrt{\frac{\log[n/\mathcal{E}(k)]}{\log(n)}}<{(1+o(1))}{M_n}\Biggr\}.
\]

Using \eqref{eq.26}, from \eqref{eq.29} we get that
\begin{equation}\label{eq.30}
\begin{aligned}[b]
\Sigma^{\rm r}_n(\xi)&\stackrel{\mathbf{P}}{=}(1+o(1))\sum_{k=1}^{\tau}
\frac{1}{\mathcal{E}(k) \sqrt{\pi \log(n/\mathcal{E}(k))}} \stackrel{\mathbf{P}}{=}
\frac{1+o(1)}{\sqrt{\pi\log(n)}}\sum_{k=1}^{\tau} \frac{1}{\mathcal{E}(k)}\\
&\stackrel{\mathbf{P}}{=}
\frac{1+o(1)}{\sqrt{\pi\log(n)}}\sum_{k=1}^{\tau}\biggl(\frac{1}{\mathcal{E}(k)}
-\frac{1}{k}\biggr) +\frac{(1+o(1))}{\sqrt{\pi\log(n)}}\sum_{k=1}^{\tau}\frac{1}{k}.
\end{aligned}
\end{equation}

In order to continue this chain of equalities, we need to analyze the stopping moment $\tau$. We can do it using various approaches, e.g. we can use results from \cite{BG}. We can easily show that if  \eqref{eq.26} is true then for $n\to \infty$ the following representation is valid:
\[%\label{tau}
\tau\stackrel{\mathbf{P}}{=} M_n+(1+o(1))[\mathcal{E}(M_n)- M_n].
\]

From here it immediately follows that for $n\to\infty$
\begin{equation}\label{eq.31}
\sum_{k=1}^{\tau}\frac{1}{k}\stackrel{\mathbf{P}}{=}\log[\mathcal{E}(M_n)]+\gamma+o(1)
\stackrel{\mathbf{P}}{=}\log(M_n)+\gamma+o(1)
\end{equation}
and
\begin{equation}\label{eq.32}
\begin{aligned}[b]
\sum_{k=1}^{\tau}\biggl(\frac{1}{\mathcal{E}(k)}-\frac{1}{k}\biggr)
&\stackrel{\mathbf{P}}{=}\sum_{k=1}^\infty\biggl(\frac{1}{\mathcal{E}(k)}
-\frac{1}{k}\biggr)
-\sum_{\tau+1}^\infty\biggl(\frac{1}{\mathcal{E}(k)}-\frac{1}{k}\biggr)\\
&\stackrel{\mathbf{P}}{=}\sum_{k=1}^\infty\biggl(\frac{1}{\mathcal{E}(k)}
-\frac{1}{k}\biggr)-\sum_{k= \tau+1}^\infty\frac{1}{k\mathcal{E}(k)}
\sum_{s=1}^k(1-e_s)\\ &\stackrel{\mathbf{P}}{=}
\sum_{k=1}^\infty\biggl(\frac{1}{\mathcal{E}(k)}-\frac{1}{k}\biggr)+ \sum_{k=
\tau+1}^\infty \frac{O(\sqrt{k})}{k^2} \stackrel{\mathbf{P}}{=}
\sum_{k=1}^\infty\biggl(\frac{1}{\mathcal{E}(k)}-\frac{1}{k}\biggr)+o(1).
\end{aligned}
\end{equation}

Thus from \eqref{eq.24}, \eqref{eq.27}, \eqref{eq.28} and \eqref{eq.30}--\eqref{eq.32}, we get that in order to prove the theorem we need to check that
\begin{equation}\label{eq.33}
\log(M_n)\biggl[\frac{1}{\sqrt{\log(b_n)}}-\frac{1}{\sqrt{\log(n)}}\biggr]
=o\biggl(\frac{1}{\sqrt{\log(n)}}\biggr).
\end{equation}
Using Taylor formula we obtain that
\[
\frac{1}{\sqrt{\log(b_n)}}-\frac{1}{\sqrt{\log(n)}}
=O\biggr(\frac{\log\sqrt{\pi\log(n)}}{\log^{3/2}(n)}\biggr).
\]
Thus \eqref{eq.33} is fulfilled if $M_n$ satisfies condition \eqref{eq.26}.\qed

\smallskip
\proofof{Theorem~\ref{main2}} Let us prove this theorem for $p=2$. We set $\lambda_1=\lambda$, and $\lambda_2=1-\lambda$. We represent $\Sigma_n(\xi)$ as follows:
\begin{equation}\label{eq.34}
\begin{aligned}[b]
\Sigma_n(\xi)&=\lambda \times \frac{1}{\lambda n}\sum_{i=1}^{\lambda n}
\exp\biggl(\frac{\xi_i^2}{2}\biggr)+ (1-\lambda)\times \frac{1}{(1-\lambda)
n}\sum_{i=\lambda n+1}^ n \exp\biggl(\frac{\xi_i^2}{2}\biggr)\\ &=\lambda
\Sigma_{\lambda n}(\xi)+(1-\lambda)\Sigma_{(1-\lambda)n}(\xi');
\end{aligned}
\end{equation}
Here, obviously $\Sigma_{\lambda n}(\xi)$ and $\Sigma_{(1-\lambda)n}(\xi')$ are independent random variables.

Using Taylor formula we can easily check that for $n\to\infty$
\[
\begin{gathered}
b_{\lambda n} =\biggl[2\log\frac{\lambda n}{\sqrt{\pi\log(\lambda n)}}\biggr]^{1/2}=
b_n +\frac{\log(\lambda)}{b_n}+ o\biggl(\frac{1}{\sqrt{\log(n)}}\biggr),\\
b_{(1-\lambda) n} =\biggl[2\log\frac{(1-\lambda) n}{\sqrt{\pi\log[(1-\lambda)
n]}}\biggr]^{1/2}= b_n +\frac{\log(1-\lambda)}{b_n}+
o\biggl(\frac{1}{\sqrt{\log(n)}}\biggr).
\end{gathered}
\]
Thus from theorem \ref{main} we get that
\[
\begin{aligned}
\lambda \Sigma_{\lambda n}(\xi) &\stackrel{\mathbf{P}}{=}\lambda
\sqrt{\frac{2}{\pi}}\biggl(b_{\lambda n} +\frac{\zeta^\circ_1}{b_{\lambda n}}\biggr)
+o\biggl(\frac{1}{\sqrt{\log(n)}}\biggr)\\
&\stackrel{\mathbf{P}}{=}\sqrt{\frac{2}{\pi}}\biggl[ \lambda b_n +\frac{\lambda
\zeta^\circ_1+\lambda\log(\lambda)}{b_n}\biggr]
+o\biggl(\frac{1}{\sqrt{\log(n)}}\biggr)
\end{aligned}
\]
and analogously
\[
\begin{aligned}
(1-\lambda) \Sigma_{(1-\lambda)n}(\xi') &\stackrel{\mathbf{P}}{=}(1-\lambda)
\sqrt{\frac{2}{\pi}}\biggl[ b_{(1-\lambda) n} +\frac{\zeta^\circ_2}{b_{(1-\lambda)
n}}\biggr] +o\biggl(\frac{1}{\sqrt{\log(n)}}\biggr) \\
&\stackrel{\mathbf{P}}{=}\sqrt{\frac{2}{\pi}}\biggl[(1-\lambda) b_n
+\frac{(1-\lambda) \zeta^\circ_2+(1-\lambda)\log(1-\lambda)}{b_n}\biggr]
+o\biggl(\frac{1}{\sqrt{\log(n)}}\biggr).
\end{aligned}
\]
This relation and \eqref{eq.34} in an obvious way completes proof of the theorem for $p=2$. The case $p>2$ can be considered analogously.\qed

\smallskip
\proofof{Theorem~\ref{lemma-b}} Let us denote for brevity
\[
\Sigma^{(-j)}(\xi)=\sum_{k\ne j}^\infty
\bar{\pi}_k^n\exp\biggl(\frac{\xi_k^2}{2}\biggr).
\]

Due to the definition of prior probabilities $\bar{\pi}_k$ we get that
\[
\Sigma^{(-j)}(\xi)\stackrel{\mathbf{P}}{=}\sum_{k=1}^\infty
\bar{\pi}_k\exp\biggl(\frac{\xi_k^2}{2}\biggr)+O\Bigl(\frac{1}{n}\Bigr),
\]
where $\bar{\pi}_k^{(-j)}$ are defined in \eqref{eq.23}. Thus from theorem \ref{main3} and Taylor formula we get the following asymptotic decomposition:
\begin{equation}\label{eq.35}
\Sigma^{(-j)}(\xi)\stackrel{\mathbf{P}}{=} \sqrt{\frac{2}{\pi}}\biggl[ b_n
+\frac{\zeta^\circ+H(\bar{\pi})}{b_n}\biggr]
+o\biggl(\frac{1}{\sqrt{\log(n)}}\biggr).
\end{equation}

It is clear that when calculating the error of the second kind without loss of generality we can assume that all $S_j$ are strictly positive and take on maximum values. More precisely, we assume that
\[
\frac{S_j}{\sigma}=r_j(n),
\]
where $r_j(n)>0$ is defined as
\[
\sqrt{\pi\log(n)}\bar{\pi}_j^n\exp\biggl[\frac{r_j^2(n)}{2}\biggr]=1,
\]
or, which is equivalent to
\begin{equation}\label{eq.36}
r_j(n)=\biggl[2\log\frac{1}{\bar{\pi}_j^n\sqrt{\pi\log(n)}}\biggr]^{1/2}.
\end{equation}

Then for conditional error of the second kind using \eqref{eq.15} and \eqref{eq.35} we obtain that
\begin{equation}\label{eq.37}
\begin{aligned}[b]
\beta_{\phi^\circ}(S\mmid\nu=j)&=\lim_{n\to\infty} \P\Biggl\{\sum_{k=1}^\infty
\bar{\pi}_k^n\exp\biggl(\frac{Y_i^2}{2\sigma^2}\biggr)\le t_\alpha^\circ(\bar{\pi}^n)
\,\Big|\,\nu=j\Biggr\}\\ &= \lim_{n\to\infty}\P\biggl\{\Sigma^{(-j)}(\xi)+
\bar{\pi}_j^n \exp\biggl(\frac{Y_j^2}{2\sigma^2}\biggr)\le
t_\alpha^\circ(\bar{\pi}^n) \biggr\}\\ &\ge\lim_{n\to\infty} \P\biggl\{\zeta^\circ\le
t_\alpha^\circ -\bar{\pi}_j^n
\sqrt{\pi\log(n)}\exp\biggl[\frac{(r_j(n)+\xi_j)^2}{2}\biggr] \biggr\}.
\end{aligned}
\end{equation}

Then we select some small value $a\in (0,1]$ and sufficiently big value  $A>1$ and continue \eqref{eq.37} as follows:
\begin{equation}\label{eq.38}
\begin{aligned}[b]
\beta_{\phi^\circ}(S\mmid \nu=j)&\ge\lim_{n\to\infty}\P\biggl\{\zeta^\circ\le
t_\alpha^\circ -\bar{\pi}_j^n
\sqrt{\pi\log(n)}\exp\biggl[\frac{(r_j(n)+\xi_j)^2}{2}\biggr]\cap -A\le\xi_j<-a
\biggr\}\\ &\ge\lim_{n\to\infty}\P\biggl\{\zeta^\circ\le t_\alpha^\circ
-\exp\biggl[\frac{A^2}{2}-r_j(n) a\biggr] \cap -A< \xi_j<-a \biggr\} \\
&\ge\P\bigl\{\zeta^\circ\le t_\alpha^\circ \bigr\} \P\bigl\{-A< \xi_j<-a
\bigr\}=(1-\alpha) \P\bigl\{-A< \xi_j<-a \bigr\}.
\end{aligned}
\end{equation}
When deriving this inequality we use the fact that $\lim\limits_{n\to\infty}r_j(n)=\infty$, which is a direct consequence of equality \eqref{eq.36}.

Since positive values $h$ and $H$ are arbitrary, inequality \eqref{eq.38} in an obvious manner completes the proof of the theorem.\qed

\Rubrika{Large Systems}

%\input{alebed.tex}

%\input{zhukov.tex}

%\Rubrika{Source Coding}

%\Rubrika{Errata}

%\Rubrika{Communication Network Theory}

%\Rubrika{Information Protection}

%\Rubrika{The International Dobrushin Prize}

%\Rubrika{Information Theory and Coding Theory}

%\Rubrika{Automata Theory}

%\Rubrika{Image Recognition}

%\Rubrika{Letters to the Editor}

%\Rubrika{Chronicles}

%\Rubrika{Obituary}

%\Rubrika{}
%
%\input{auth-ind.tex}
%
%\input{contents.tex}

\end{document}